\documentclass[11pt]{amsart}
\usepackage{amsfonts}
\usepackage{amsmath,amsthm,amsfonts,amssymb}
\newtheorem{lemma}{Lemma}
\newtheorem{theorem}{Theorem}
\newtheorem{corollary}{Corollary}

\def\bl{\begin{lemma}}
\def\bt{\begin{theorem}}
\def\el{\end{lemma}}
\def\et{\end{theorem}}
\def\bp{\begin{proof}}
\def\ep{\end{proof}}
\def\bc{\begin{corollary}}
\def\ec{\end{corollary}}
\def\iy{\infty}

\def\mc{\mathcal}

\def\su{\subset}
\def\l{\lambda}

\def\o{\omega}
\def\O{\Omega}
\def\a{\alpha}
\def\b{\beta}
\def\p{\partial}
\def\su{\subset}
\def\G{\Gamma}
\def\s{\sigma}
\def\d{\delta}

\def\r{\rho}

\def\-{\setminus}
\def\s{\sigma}

\def\ov{\overline}
\def\lt{\left}
\def\rt{\right}
\def\+{\bigcup}
\def\.{\bigcap}

\title[A note on Schwarz-Pick lemma]{A note on Schwarz-Pick lemma for bounded complex-valued
harmonic functions in the unit ball of $\mathbb{R}^n$}
\thanks{Research supported by the National Natural Science Foundation of China
(No. 11201199) and by the Scientific Research Foundation of Jinling
Institute of Technology (No. Jit-b-201221).}

\author {Shaoyu Dai and Yifei Pan}

\address{Department of General Study Program, Jinling Institute of
Technology, Nanjing 211169, China}
\address{\it E-mail address: dymdsy@163.com}

\address{Department of Mathematical Sciences, Indiana University -
Purdue University Fort Wayne, Fort Wayne, IN 46805-1499, USA}
\address{\it E-mail address: pan@ipfw.edu}

\begin{document}

\numberwithin{equation}{section}

\begin{abstract}
In this paper we prove a Schwarz-Pick lemma for bounded
complex-valued harmonic functions in the unit ball of
$\mathbb{R}^n$.
\end{abstract}

\maketitle



\section{Introduction}

This paper is a note about Chen's paper \cite{Chen}. Using the same
method in \cite{Chen}, we obtain a result Theorem 1, which extends
the Schwarz-Pick lemma \cite{Chen} for planar harmonic mappings to
bounded complex-valued harmonic functions in the unit ball of
$\mathbb{R}^n$. In addition, motivated by \cite{Chen} and this
paper, we consider a Schwarz lemma for harmonic mappings between
real unit balls in another paper. Now we introduce some denotation
and the background.

Let $n$ be a positive integer greater than 1. $\mathbb{R}^n$ is the
real space of dimension $n$. For
$x=(x_1,\cdots,x_n)\in\mathbb{R}^n$, let
$|x|=(|x_1|^2+\cdots+|x_n|^2)^{1/2}$. Let
$\mathbb{B}^n=\{x\in\mathbb{R}^n: |x|<1\}$ be the unit ball of
$\mathbb{R}^n$. The unit sphere, the boundary of $\mathbb{B}^n$ is
denoted by $S$; normalized surface-area measure on $S$ is denoted by
$\s$ (so that $\s(S)=1$). Let $S^+$ denote the northern hemisphere
$\{x=(x_1,\cdots,x_n)\in S: x_n>0\}$ and let $S^-$ denote the
southern hemisphere $\{x=(x_1,\cdots,x_n)\in S: x_n<0\}$.
$N=(0,\cdots,0,1)$ denotes the north pole of $S$.
$B_r=\{x\in\mathbb{R}^n: |x|<r\}$ is the open ball centered at
origin of radius $r$; its closure is the closed ball $\ov B_r$. A
twice continuously differentiable, complex-valued function $F$
defined on $\mathbb{B}^n$ is harmonic on $\mathbb{B}^n$ if and only
if $\Delta F\equiv0$, where $\Delta=D_1^2+\cdots+D_n^2$ and $D_j^2$
denotes the second partial derivative with respect to the $j^{th}$
coordinate variable $x_j$. By $\O_n$, we denote the class of all
complex-valued harmonic functions $F(x)$ on $\mathbb{B}^n$ with
$|F(x)|<1$ for $x\in\mathbb{B}^n$.

Let $\mathbb{D}$ be the unit disk in the complex plane $\mathbb{C}$.
Denote the disk $\{z\in\mathbb{C}: |z<r|\}$ by $D_r$; its closure is
the closed disk $\ov D_r$.

For a holomorphic function $f$ from $\mathbb{D}$ into $\mathbb{D}$,
the classical Schwarz lemma says that if $f(0)=0$, then
\begin{equation}\label{Slemma1}
|f(z)|\leq|z|
\end{equation}
holds for $z\in\mathbb{D}$. For $0<r<1$, \eqref{Slemma1} may be
written in the following form:
\begin{equation}\label{Slemma2}
f(\ov D_r)\subset\ov D_r.
\end{equation}
So the classical Schwarz lemma can be regarded as considering the
region of $f(\ov D_r)$. If the condition  $f(0)=0$ is relaxed, then
what the region of $f(\ov D_r)$ is. The answer can be found in the
classical Schwarz-Pick lemma. By Schwarz-Pick lemma \cite{Ahlfors},
it is known that
\begin{equation}\label{SPlemma1}
\frac{|f(z_1)-f(z_2)|}{|1-\ov{f(z_2)}f(z_1)|}\le\frac{|z_1-z_2|}{|1-\ov{z_2}z_1|}
\end{equation}
holds for $z_1,z_2\in\mathbb{D}$. Using the notations
$$d_p(z_1,z_2)=\frac{|z_1-z_2|}{|1-\ov{z_2}z_1|}$$
for the pseudo-distance between $z_1,z_2\in\mathbb{D}$, we know that
\begin{equation}\label{SPlemma2}
d_p(f(z_1),f(z_2))\leq d_p(z_1,z_2)
\end{equation}
for $z_1,z_2\in\mathbb{D}$ by \eqref{SPlemma1}. Denote
$\ov{\triangle}(z,r)=\{\zeta\in\mathbb{D}: d_p(\zeta,z)\leq r,
z\in\mathbb{D}, 0<r<1$ for the closed pseudo-disk with center at $z$
and pseudo-radius $r$. Then \eqref{SPlemma2} may be written in the
following form:
$$f(\ov{\triangle}(z,r))\subset\ov{\triangle}(f(z),r)$$
for $z\in\mathbb{D}$ and $0<r<1$. Note that $\ov{\triangle}(0,r)=\ov
D_r$. So for $f$ without the assumption $f(0)=0$, we know
\begin{equation}\label{SPlemma3}
f(\ov D_r)\subset\ov{\triangle}(f(0),r).
\end{equation}
When $f(0)=0$, \eqref{SPlemma3} becomes \eqref{Slemma2}.

For a complex-valued harmonic function $F$ on $\mathbb{D}$ such that
$F(\mathbb{D})\subset\mathbb{D}$ and $F(0)=0$, it is known
\cite{Heinz} that
\begin{equation}\label{HSlemma1}
F(z)\leq\frac{4}{\pi}\arctan|z|
\end{equation}
holds for $z\in\mathbb{D}$. For $0<r<1$, \eqref{HSlemma1} may be
written in the following form:
\begin{equation}\label{HSlemma2}
F(\ov D_r)\subset\ov D_{\frac{4}{\pi}\arctan r}.
\end{equation}
If the condition  $F(0)=0$ is relaxed, then what the region of
$F(\ov D_r)$ is. Unfortunately, the composition $f\circ F$ of a
harmonic function $F$ and a holomorphic function $f$ do not need to
be harmonic, so it is a serious problem to seek the estimate
corresponding to \eqref{SPlemma3} for a harmonic function $F$
without the assumption $F(0)=0$. Fortunately, Chen resolved this
problem in \cite{Chen}. In \cite{Chen}, for any $0<r<1$ and
$0\leq\rho<1$, the author construct a closed domain $E_{r,\r}$,
which contains $\r$ and is symmetric to the real axis, with the
following properties: Let $z\in\mathbb{D}$ and $w=\r e^{i\a}$ be
given. For every complex-valued harmonic function $F$ with
$F(\mathbb{D})\subset\mathbb{D}$ and $F(z)=w$, the author has
$F(\ov{\triangle}(z,r))\subset e^{i\a}E_{r,\r}=\{e^{i\a}\zeta:
\zeta\in E_{r,\r}\}$; conversely, for every $w'\in e^{i\a}E_{r,\r}$,
there exists a complex-valued harmonic function $F$ such that
$F(\mathbb{D})\subset\mathbb{D}$, $F(z)=w$ and $F(z')=w'$ for some
$z'\in\partial\triangle(z,r)$. Obviously, by Chen's result, we know
that for a complex-valued harmonic function $F$ on $\mathbb{D}$ such
that $F(\mathbb{D})\subset\mathbb{D}$ without the assumption
$F(0)=0$, if $F(0)=\r e^{i\a}$, then
\begin{equation}\label{Chen}
F(\ov D_r)\subset e^{i\a}E_{r,\r},\
\end{equation}
which is sharp. \eqref{Chen} is the estimate for complex-valued
harmonic functions corresponding to \eqref{SPlemma3}. Note that a
complex-valued harmonic function $F$ on $\mathbb{D}$ such that
$F(\mathbb{D})\subset\mathbb{D}$ can be seen as $F\in\O_2$. So it is
natural to consider the same problem in $\O_n$.

For $F\in\O_n$, harmonic Schwarz lemma \cite{Axler} says that if
$F(0)=0$, then
\begin{equation}\label{HSlemma3}
|F(x)|\leq U(|x|N)
\end{equation}
holds for $x\in\mathbb{B}^n$, where $U$ is the Poisson integral of
the function that equals 1 on $S^+$ and -1 on $S^-$. For $0<r<1$,
\eqref{HSlemma3} may be written in the following form:
\begin{equation}\label{HSlemma4}
F(\ov B_r)\subset\ov D_{U(rN)}.
\end{equation}
If the condition  $F(0)=0$ is relaxed, then what the region of
$F(\ov B_r)$ is. This problem will be solved in this paper.

In this paper, by the same method in \cite{Chen}, we obtain the
following theorem about the region of $F(\ov B_r)$. The result is
sharp. When $n=2$, our result is coincident with \eqref{Chen}. And
when $F(0)=0$, our result is coincident with \eqref{HSlemma4}. Note
that in the following theorem, $E_{r,\r}$ is defined as \eqref{Err}.

\bt Let $0\leq\r<1$, $\a\in\mathbb{R}$ and $0<r<1$ be given. Then
for every harmonic function $F$ with
$F(\mathbb{B}^n)\subset\mathbb{D}$ and $F(0)=\r e^{i\a}$, we have
$F(\ov B_r)\su e^{i\a}E_{r,\r}=\{e^{i\a}\zeta: \zeta\in E_{r,\r}\}$;
conversely, for every  $w'\in e^{i\a}E_{r,\r}$, there exists a
harmonic function $F$ such that $F(\mathbb{B}^n)\subset\mathbb{D}$,
$F(0)=\r e^{i\a}$ and $F(rN)=w'$.\et

The theorem above will be proved by three steps as follows.\\
Step 1: find the extremal line of $F(\ov B_r)$ in the normal
direction of $e^{0i}$, which is related to the value of $F(0)$.\\
Step 2: find the extremal line of $F(\ov B_r)$ in the normal
direction of a given direction. For a given direction of $e^{i\b}$
with $-\pi\le\b\le\pi$, construct a new harmonic function
$F_{\b}=e^{-i\b}F$ through rotating $F(\ov B_r)$ by an
anti-clockwise rotation of angle $\b$. Using the result of Step1, we
will have the the extremal line of $F_{\b}(\ov B_r)$ in the normal
direction of $e^{0i}$, which is denoted by $l'_\b$. Note that $F(\ov
B_r)$ can be obtained from $F_{\b}(\ov B_r)$ by a clockwise rotation
of angle $\b$. Then the extremal line of $F(\ov B_r)$ in the normal
direction of $e^{i\b}$, which is denoted by $l_\b$, can be obtained
from $l'_\b$ by a clockwise rotation of angle $\b$.\\
Step 3: using the result of Step2, we will obtain all the extremal
lines of $F(\ov B_r)$ in every normal direction, with which we can
wrap $F(\ov B_r)$ and obtain the region of $F(\ov B_r)$.\\
Step 1 will be solved in Section 2. Step 2 and Step 3 will be solved
in Section 3.

\section{Some lemmas}

In this section, we will introduce some lemmas, which are important
for the proof of Theorem 2. Lemma 1 will be used in Lemma 2. Lemma 2
will be used in Lemma 3. Lemma 3 and Lemma 4 will be used in Theorem
2.

Now we give Lemma 1 first.  Lemma 1 constructs a bijection $(R,I)$
from $\mathbb{R}\times\mathbb{R}^+$ onto the upper half disk
$\{(a,b):a\in\mathbb{R},b\in\mathbb{R},a^2+b^2<1,\ b>0\}$, which
will be used to construct $u_{a,b,r}$ in Lemma 2 for the case that
$b>0$.

For $0<r<1$, $\mu>0$ and real number $\l$, define
\begin{equation}\label{Ady}
A_{r,\l,\mu}(\omega) =\frac1\mu\lt(\frac{1}{|rN-\o|^n}-\l\rt),\ \ \
\ \ \ \o\in S,
\end{equation}
and
\begin{equation}\label{RIdy}
R(r,\l,\mu)=\int_S\frac{A_{r,\l,\mu}(\o)}
{\sqrt{1+A^2_{r,\l,\mu}(\o)}}\,d\s, \quad I(r,\l,\mu)=\int_S\frac{1}
{\sqrt{1+A^2_{r,\l,\mu}(\o)}}\,d\s.
\end{equation}
The idea of the conformation of $A_{r,\l,\mu}(\omega)$,
$R(r,\l,\mu)$ and $I(r,\l,\mu)$ originates from the needs of
\eqref{u0dy} and \eqref{u4}.

\bl\label{lemmaRI} Let $0<r<1$ be fixed. Then, there exist a unique
pair of real functions $\l=\l(r,a,b)$ and $\mu=\mu(r,a,b)>0$,
defined on the upper half disk $\{(a,b):a^2+b^2<1,\ b>0\}$ and
analytic in the real sense, such that $R(r,\l(r,a,b),\mu(r,a,b))=a$
and $I(r,\l(r,a,b),\mu(r,a,b))=b$ for any point $(a,b)$ in the half
disk. \el

\bp A simple calculation gives
\begin{equation}\label{1}
\frac{\p R(r,\l,\mu)}{\p\l}=-\frac1{\mu}\int_S
\frac1{(1+A^2_{r,\l,\mu}(\o))^{3/2}}d\s,
\end{equation}
\begin{equation}\label{2}
\frac{\p R(r,\l,\mu)}{\p\mu}=-\frac1{\mu}\int_S
\frac{A_{r,\l,\mu}(\o)}{(1+A^2_{r,\l,\mu}(\o))^{3/2}}d\s,
\end{equation}
\begin{equation}\label{3}
\frac{\p I(r,\l,\mu)}{\p\l}=\frac1{\mu}\int_S
\frac{A_{r,\l,\mu}(\o)}{(1+A^2_{r,\l,\mu}(\o))^{3/2}}d\s,
\end{equation}
\begin{equation}\label{4}
\frac{\p I(r,\l,\mu)}{\p\mu}=\frac1{\mu}\int_S
\frac{A^2_{r,\l,\mu}(\o)}{(1+A^2_{r,\l,\mu}(\o))^{3/2}}d\s.
\end{equation}
It is easy to see that\\
(i) by \eqref{1}, $\p R(r,\l,\mu)/\p\l<0$ for any $\l$ and $\mu>0$,
and $R(r,\l,\mu)$ is strictly decreasing as a function of $\l$ for a
fixed $\mu$;\\
(ii) by \eqref{RIdy}, for a fixed $\mu$, $R(r,\l,\mu)\to-1$ or $1$
according to
$\l\to+\iy$ or $\l\to-\iy$;\\
(iii) by \eqref{1} - \eqref{4} and the convexity of the square
function, $$\frac{\p R(r,\l,\mu)}{\p\l}\frac{\p
I(r,\l,\mu)}{\p\mu}-\frac{\p R(r,\l,\mu)}{\p\mu}\frac{\p
I(r,\l,\mu)}{\p\l}<0$$ for any $\l$ and
$\mu>0$;\\
(iiii) by \eqref{RIdy}, $0<I(r,\l,\mu)<1$ for any $\l$ and $\mu>0$.

By (i) and (ii), we know that for fixed $\mu$, $R_(r,\l,\mu)$ is
strictly decreasing from $1$ to $-1$ as $\l$ increasing from $-\iy$
to $+\iy$. Then for any $-1<a<1$ and fixed $\mu$, there exists a
unique real number $\l(\mu,a)$ such that
\begin{equation}\label{wyy}
R(r,\l,\mu)\lt|_{\l=\l(\mu,a)}\rt.=a.
\end{equation}
Further, using the implicit function theorem, we have that the
function $\l=\l(\mu,a)$ defined on $\{(\mu,a):\mu>0,-1<a<1\}$ is a
continuous function and
$$\frac{\p\l(\mu,a)}{\p\mu}=-\lt.\lt(\frac{\p R(r,\l,\mu)}{\p\mu}
\lt/\frac{\p R(r,\l,\mu)}{\p\l}\rt.\rt)\rt|_{\l=\l(\mu,a)}.$$

Next, we consider the function $I(r,\l(\mu,a),\mu)$ for $\mu>0$.
\begin{equation*}
\begin{split}
&\ \ \frac{\p I(r,\l(\mu,a),\mu)}{\p\mu}\\
&=\lt.\lt(\frac{\p
I(r,\l,\mu)}{\p\l}\frac{\p\l(\mu,a)}{\p\mu}+\frac{\p
I(r,\l,\mu)}{\p\mu}\rt)\rt|_{\l=\l(\mu,a)}\\
&=\lt.\lt(\lt(\frac{\p R(r,\l,\mu)}{\p\l}\frac{\p
I(r,\l,\mu)}{\p\mu}-\frac{\p R(r,\l,\mu)}{\p\mu}\frac{\p
I(r,\l,\mu)}{\p\l}\rt)\lt/\frac{\p
R(r,\l,\mu)}{\p\l}\rt.\rt)\rt|_{\l=\l(\mu,a)}.
\end{split}
\end{equation*}
By (i) and (iii), we have $\frac{\p I(r,\l(\mu,a),\mu)}{\p\mu}>0$,
which shows that $I(r,\l(\mu,a),\mu)$ is strictly increasing as a
function of $\mu$ on $(0,+\infty)$ for a fixed $a$. Note that
(iiii). Thus, for a fixed $a$, $I(r,\l(\mu,a),\mu)$ respectively has
finite limit as $\mu\to 0$ and as $\mu\to+\iy$.

For a fixed $a$, we claim that $I(r,\l(\mu,a),\mu)\to 0$ as $\mu\to
0$, and $I(r,\l(\mu,a),\mu)\to\sqrt{1-a^2}$ as $\mu\to+\iy$.

As $\mu\to0$, there exists a subsequence $\mu_k\to0$ such that
$\l(\mu_k,a)$ has a finite limit $t$ or tend to $\iy$. We only need
to prove that $I(r,\l(\mu_k,a),\mu_k)\to0$ as $k\to\iy$. Since
$I(r,\l(\mu_k,a),\mu_k)=\int_S\frac{1}
{\sqrt{1+A^2_{r,\l(\mu_k,a),\mu_k}(\o)}}\,d\s$, we only need to
prove that $|A_{r,\l(\mu_k,a),\mu_k}(\o)|\to+\iy$ almost everywhere
on $S$. Note that
$$|A_{r,\l(\mu_k,a),\mu_k}(\o)|=\frac{1}{\mu_k}\lt|\frac{1}{|rN-\o|^n}-\l(\mu_k,a)\rt|$$ and
$$\frac{1}{(1+r)^n}\leq\frac{1}{|rN-\o|^n}\leq\frac{1}{(1-r)^n}.$$
If $\l(\mu_k,a)\to t$ as $k\to\iy$, then
$\frac{1}{|rN-\o|^n}-\l(\mu_k,a)$ is bounded and
$\frac{1}{|rN-\o|^n}-\l(\mu_k,a)\neq0$ almost everywhere on $S$.
Thus $|A_{r,\l(\mu_k,a),\mu_k}(\o)|\to+\iy$ almost everywhere on
$S$. If $\l(\mu_k,a)\to\iy$ as $k\to\iy$, then it is obvious that
$|A_{r,\l(\mu_k,a),\mu_k}(\o)|\to+\iy$ uniformly for $\o\in S$. The
first claim is proved.

As $\mu\to+\iy$, $\frac{1}{\mu}\frac{1}{|rN-\o|^n}\to0$ uniformly
for $\o\in S$. If there exists a subsequence $\mu_k\to+\iy$ such
that $\l(\mu_k,a)/\mu_k\to\iy$, then
$|A_{r,\l(\mu_k,a),\mu_k}(\o)|\to+\iy$ uniformly for $\o\in S$, and
$I(r,\l(\mu_k,a),\mu_k)\to 0$, a contradiction. This shows that for
$\l(\mu,a)/\mu$ is bounded as $\mu\to+\iy$. Thus there exists a
subsequence $\mu_k\to+\iy$ such that $-\l(\mu_k,a)/\mu_k$ tend to a
finite limit $t$. That is
\begin{equation}\label{lmu}
\lim_{k\to\iy}-\l(\mu_k,a)/\mu_k=t.
\end{equation}
we only need to prove that $I(r,\l(\mu_k,a),\mu_k)\to\sqrt{1-a^2}$
as $k\to\iy$. Let $(A(\o))_k=A_{r,\l(\mu_k,a),\mu_k}(\o)$. By
\eqref{Ady}, \eqref{lmu} and $\mu_k\to+\iy$, we obtain
\begin{equation}\label{wR}
\begin{split}
\lim_{k\to\iy}\frac{(A(\o))_k}{\sqrt{1+((A(\o))_k)^2}}
&=\lim_{k\to\iy}\frac{\frac1\mu_k\lt(\frac{1}{|rN-\o|^n}-\l(\mu_k,a)\rt)}
{\sqrt{1+\frac1{\mu_k^2}\lt(\frac{1}{|rN-\o|^n}-\l(\mu_k,a)\rt)^2}}\\
&=\lim_{k\to\iy}\frac{-\frac{\l(\mu_k,a)}{\mu_k}}
{\sqrt{1+\lt(\frac{\l(\mu_k,a)}{\mu_k}\rt)^2}}
=\frac{t}{\sqrt{1+t^2}}
\end{split}
\end{equation}
uniformly for $\o\in S$, and
\begin{equation}\label{wI}
\begin{split}
\lim_{k\to\iy}\frac{1}{\sqrt{1+((A(\o))_k)^2}}
&=\lim_{k\to\iy}\frac{1}
{\sqrt{1+\frac1{\mu_k^2}\lt(\frac{1}{|rN-\o|^n}-\l(\mu_k,a)\rt)^2}}\\
&=\lim_{k\to\iy}\frac{1}
{\sqrt{1+\lt(\frac{\l(\mu_k,a)}{\mu_k}\rt)^2}}
=\frac{1}{\sqrt{1+t^2}}
\end{split}
\end{equation}
uniformly for $\o\in S$. By the Lebesgue's dominated convergence
theorem and \eqref{RIdy}, \eqref{wR}, \eqref{wI} we have
\begin{equation}\label{wjR}
\begin{split}
\lim_{k\to\iy}R(r,\l(\mu_k,a),\mu_k)
&=\lim_{k\to\iy}\int_S\frac{(A(\o))_k}{\sqrt{1+((A(\o))_k)^2}}\,d\s\\
&=\int_S\lim_{k\to\iy}\frac{(A(\o))_k}{\sqrt{1+((A(\o))_k)^2}}\,d\s\\
&=\frac{t}{\sqrt{1+t^2}},
\end{split}
\end{equation}
and
\begin{equation}\label{wjI}
\begin{split}
\lim_{k\to\iy}I(r,\l(\mu_k,a),\mu_k)
&=\lim_{k\to\iy}\int_S\frac{1}{\sqrt{1+((A(\o))_k)^2}}\,d\s\\
&=\int_S\lim_{k\to\iy}\frac{1}{\sqrt{1+((A(\o))_k)^2}}\,d\s\\
&=\frac{1}{\sqrt{1+t^2}}.
\end{split}
\end{equation}
Note that $R(r,\l(\mu_k,a),\mu_k)\equiv a$ by \eqref{wyy}, and
$\lt(\frac{t}{\sqrt{1+t^2}}\rt)^2
+\lt(\frac{1}{\sqrt{1+t^2}}\rt)^2=1$. Then by \eqref{wjR} we obtain
that $\frac{t}{\sqrt{1+t^2}}=a$ and
$\frac{1}{\sqrt{1+t^2}}=\sqrt{1-a^2}$. Consequently by \eqref{wjI},
$$\lim_{k\to\iy}I(r,\l(\mu_k,a),\mu_k)=\sqrt{1-a^2}.$$
The second claim is proved.

It is proved that $I(r,\l(\mu,a),\mu)$ is continuous and strictly
increasing from $0$ to $\sqrt{1-a^2}$ as $\mu$ increasing from 0 to
$+\iy$. Thus, for any $0<b<\sqrt{1-a^2}$ and $-1<a<1$, there exists
a unique real number $\mu(a,b)$ such that
\begin{equation}\label{muab}
I(r,\l(\mu(a,b),a),\mu(a,b))=b.
\end{equation}
Further, using the implicit function theorem, we have the function
$\mu(a,b)$ defined on $\{(a,b):a^2+b^2<1,\ b>0\}$ is a continuous
function.

Denote $\l(\mu(a,b),a)$ by $\l(r,a,b)$. Denote $\mu(a,b)$ by
$\mu(r,a,b)$. We have proved that there exist a unique pair of
functions $\l=\l(r,a,b)$ and $\mu=\mu(r,a,b)$ such that
$$R(r,\l(r,a,b),\mu(r,a,b))=a, \ \ I(r,\l(r,a,b),\mu(r,a,b))=b$$
on the upper half disk. The real analyticity of $\l=\l(r,a,b)$ and
$\mu=\mu(r,a,b)$ is asserted by the implicit function theorem. The
lemma is proved. \ep

Let $a$ and $b$ be two numbers such that $0\le b<1$, $-1<a<1$ and
$a^2+b^2<1$. Let ${\mc U}_{a,b}$ denote the class of real-valued
functions $u\in L^\infty(S)$ satisfying the following conditions:
\begin{equation} \label{uyq}
\|u\|_{\infty}\le 1,\quad \int_S u(\o)d\s=a,\quad
\int_S\sqrt{1-u^2(\o)}d\s\ge b.
\end{equation}
Every function $u\in L^\infty(S)$ defines a harmonic function
$$U(x)=\int_S\frac{1-|x|^2}{|x-\o|^n}
u(\o)d\s\quad\mbox{for}\quad x\in\mathbb{B}^n.$$ Let $0<r<1$ and
define a functional $L_r$ on $L^\infty(S)$ by
\begin{equation} \label{Ldy}
L_r(u)=U(rN) =\int_S\frac{1-r^2}{|rN-\o|^n} u(\o)d\s.
\end{equation}

Obviously, ${\mc U}_{a,b}$ is a closed set, and $L_r$ is a
continuous functional on ${\mc U}_{a,b}$. Then there exists a
extremal function such that $L_r$ attains its maximum on ${\mc
U}_{a,b}$ at the extremal function. We will claim in the following
lemma that the extremal function is unique. In the proof of the
following lemma, we will construct a function $u_0$ first and then
prove that $u_0$ is the unique extremal function, which will be
denoted by $u_{a,b,r}$.

\bl\label{lemmauabr} For any $a$, $b$ and $r$ satisfying the above
conditions, there exists a unique extremal function
$u_{a,b,r}\in{\mc U}_{a,b}$ such that $L_r$ attains its maximum on
${\mc U}_{a,b}$ at $u_{a,b,r}$. \el

\bp Let $a$, $b$ and $r$ be fixed. First assume that $b>0$. From
Lemma \ref{lemmaRI}, we have $\l=\l(r,a,b)$ and $\mu=\mu(r,a,b)>0$
such that $R(r,\l,\mu)=a$ and $I(r,\l,\mu)=b.$ For the need of
\eqref{u4}, let
\begin{equation} \label{u0dy}
u_0(\o)=\frac{A_{r,\l,\mu}(\o)}{\sqrt{1+A^2_{r,\l,\mu}(\o)}},
\end{equation}
where $A_{r,\l,\mu}(\o)$ is defined as \eqref{Ady}. Then
$||u_0||_\iy<1$ and by \eqref{RIdy}, we know
\begin{equation} \label{u0td}
\int_S u_0(\o)d\s=R(r,\l,\mu)=a,\quad \int_S \sqrt{1-u_0^2(\o)}d\s
=I(r,\l,\mu)=b.
\end{equation}
This means that $u_0\in {\mc U}_{a,b}$.

Let $u\in{\mc U}_{a,b}$. By \eqref{uyq} and \eqref{u0td}, we have
\begin{equation} \label{u1}
\l\int_S (u_0(\o)-u(\o))d\s=0,
\end{equation}
\begin{equation} \label{u2}
\mu\int_S (\sqrt{1-u^2_0(\o)}-\sqrt{1-u^2(\o)})d\s\leq0.
\end{equation}
By the Taylor formula of the function $\sqrt{1-x^2}$, we have
\begin{equation} \label{u3}
\begin{split}
\sqrt{1-u^2(\o)}-\sqrt{1-u^2_0(\o)}&=
\frac{u_0(\o)(u_0(\o)-u(\o))}{\sqrt{1-u^2_0(\o)}}-\frac{(u_0(\o)-u(\o))^2}{2(1-\xi^2)^{3/2}}\\
&\leq\frac{u_0(\o)(u_0(\o)-u(\o))}{\sqrt{1-u^2_0(\o)}},
\end{split}
\end{equation}
where $\xi$ is a real number between $u_0(\o)$ and $u(\o)$. By
\eqref{u0dy} and \eqref{Ady}, we have
\begin{equation} \label{u4}
\frac{1}{|rN-\o|^n}-\l-\frac{\mu u_0(\o)}{\sqrt{1-u^2_0(\o)}}=0.
\end{equation}
Then by \eqref{Ldy} and \eqref{u1}-\eqref{u4}, we obtain that
\begin{eqnarray*}
&&\frac{L_r(u_0)-L_r(u)}{1-r^2}=\int_S\frac{u_0(\o)-u(\o)}{|rN-\o|^n}d\s\\
&\geq&\int_S\frac{u_0(\o)-u(\o)}{|rN-\o|^n}d\s-\l\int_S
(u_0(\o)-u(\o))d\s-\mu\int_S
(\sqrt{1-u^2(\o)}-\sqrt{1-u^2_0(\o)})d\s\\
&=&\int_S\frac{u_0(\o)-u(\o)}{|rN-\o|^n}d\s-\l\int_S
(u_0(\o)-u(\o))d\s-\mu\int_S
\frac{u_0(\o)(u_0(\o)-u(\o))}{\sqrt{1-u^2_0(\o)}}d\s\\
&\ &\ \ \ \ \ \ \ \ \ \ \ \ \ \ \ \ \ \ \ \ +\mu\int_S
\frac{(u_0(\o)-u(\o))^2}{2(1-\xi^2)^{3/2}}d\s\\
&\geq&\int_S\frac{u_0(\o)-u(\o)}{|rN-\o|^n}d\s-\l\int_S
(u_0(\o)-u(\o))d\s-\mu\int_S\frac{
u_0(\o)(u_0(\o)-u(\o))}{\sqrt{1-u^2_0(\o)}}d\s\\
&=&\int_S (u_0(\o)-u(\o)) \lt(\frac{1}{|rN-\o|^n}-\l-\frac{\mu
u_0(\o)}{\sqrt{1-u^2_0(\o)}}\rt)d\s\\
&=&0.
\end{eqnarray*}
Thus $L_r(u_0)\geq L_r(u)$ with equality if and only if $\mu\int_S
\frac{(u_0(\o)-u(\o))^2}{2(1-\xi^2)^{3/2}}d\s=0$. Therefore
$L_r(u_0)\geq L_r(u)$ with equality if and only if $u(\o)=u_0(\o)$
almost everywhere. This shows that $u_0(\o)$ is the unique extremal
function, which will be denoted by $u_{a,b,r}(\o)$.

Next we consider the case that $b=0$. For a real number $d$, let
\begin{equation} \label{Sd1}
S_d=\{x\in S: |N-x|=d\},
\end{equation}
\begin{equation} \label{Sd2}
S^+_d=\{x\in S: |N-x|<d\},
\end{equation}
\begin{equation} \label{Sd3}
S^-_d=\{x\in S: |N-x|>d\}.
\end{equation}
For a fixed real number $a$ such that $-1<a<1$, there exists a
unique real number $d_a$ such that $\s(S^+_{d_a})=\frac{1+a}{2}$ and
$\s(S^-_{d_a})=\frac{1-a}{2}$. Let
\begin{equation}\label{b=0u0dy}
u_0(\o)=
\begin{cases}
1, &\o\in S^+_{d_a};\\
0, &\o\in S_{d_a};\\
-1, &\o\in S^-_{d_a}.
\end{cases}
\end{equation}
We want to prove that $u_0$ is just the unique extremal function,
which will be denoted by $u_{a,0,r}(\o)$.

It is obvious that $u_0\in\mc U_{a,0}$. Let $u\in{\mc U}_{a,0}$. By
\eqref{uyq} and \eqref{b=0u0dy}, we have
\begin{equation} \label{xy1}
\int_S (u_0(\o)-u(\o))d\s=0,
\end{equation}
\begin{equation} \label{xy2}
u_0(\o)-u(\o)\geq0\quad\mbox{for}\quad\o\in S^+_{d_a},
\end{equation}
\begin{equation} \label{xy3}
u_0(\o)-u(\o)\leq0\quad\mbox{for}\quad\o\in S^-_{d_a},
\end{equation}
Let
\begin{equation} \label{xy4}
J_a=|rN-x_0|, \ \ \mbox{where $x_0\in S_{d_a}$}.
\end{equation}
Note that
\begin{equation} \label{xy5}
|rN-\o|<J_a\quad\mbox{for}\quad\o\in S^+_{d_a},
\end{equation}
\begin{equation} \label{xy6}
|rN-\o|>J_a\quad\mbox{for}\quad\o\in S^-_{d_a}.
\end{equation}
Then by \eqref{Ldy} and \eqref{xy1}-\eqref{xy6}, we obtain that
\begin{equation*}
\begin{split}
&\frac{L_r(u_0)-L_r(u)}{1-r^2}\\
&=\int_S\frac{u_0(\o)-u(\o)}{|rN-\o|^n}d\s\\
&=\int_S(\frac{1}{|rN-\o|^n}-\frac{1}{J^n_a})(u_0(\o)-u(\o))d\s\\
&=\int_{S^+_{d_a}}(\frac{1}{|rN-\o|^n}-\frac{1}{J^n_a})(u_0(\o)-u(\o))d\s
+\int_{S^-_{d_a}}(\frac{1}{|rN-\o|^n}-\frac{1}{J^n_a})(u_0(\o)-u(\o))d\s\\
&\geq0.
\end{split}
\end{equation*}
Thus $L_r(u_0)\geq L_r(u)$ with equality if and only if
$u(\o)=u_0(\o)$ almost everywhere. The lemma is proved. \ep

Let $a$ and $b$ be two real numbers with $a^2+b^2<1$, and $0<r<1$.
If $b\ge0$, $u_{a,b,r}$ has been defined in Lemma \ref{lemmauabr}.
Now, define
\begin{equation} \label{vdy}
v_{a,b,r}(\o)=\sqrt{1-u_{a,b,r}^2(\o)}\quad\mbox{for}\quad\o\in S,
\end{equation}
and
\begin{equation} \label{Udy}
U_{a,b,r}(x)=\int_S\frac{1-|x|^2}{|x-\o|^n} u_{a,b,r}(\o)d\s,
\end{equation}
\begin{equation} \label{Vdy}
V_{a,b,r}(x)=\int_S\frac{1-|x|^2}{|x-\o|^n} v_{a,b,r}(\o)d\s.
\end{equation}
For $b<0$, let
\begin{equation} \label{UVtd}
U_{a,b,r}(x)=U_{a,-b,r}(x),\quad V_{a,b,r}(x)=-V_{a,-b,r}(x).
\end{equation}
Then for any $a\in\mathbb{R}$, $b\in\mathbb{R}$ and $a^2+b^2<1$, let
\begin{equation} \label{Fabrdy}
F_{a,b,r}(x)=U_{a,b,r}(x)+iV_{a,b,r}(x)\quad\mbox{for}\quad
x\in\mathbb{B}^n.
\end{equation}
The harmonic function $F_{a,b,r}(x)=U_{a,b,r}(x)+iV_{a,b,r}(x)$
satisfies $F_{a,b,r}(0)=a+bi$ and
$F_{a,b,r}(\mathbb{B}^n)\subset\mathbb{D}$, since we will show that
$|U_{a,b,r}(x)|^2+|V_{a,b,r}(x)|^2<1$. By the convexity of the
square function,
$$|U_{a,b,r}(x)|^2+|V_{a,b,r}(x)|^2\leq\int_S\frac{1-|x|^2}{|x-\o|^n}
(u_{a,b,r}^2(\o)+v_{a,b,r}^2(\o))d\s=1$$ with equality if and only
if $u_{a,b,r}(\o)$ and $v_{a,b,r}(\o)$ are constants almost
everywhere on $S$. However $u_{a,b,r}(\o)$ is not possible a
constant almost everywhere on $S$. Thus
$|U_{a,b,r}(x)|^2+|V_{a,b,r}(x)|^2<1$.

The functions $F_{a,b,r}$ are the extremal functions in the
following lemma.

\bl\label{lemma3} Let $F(x)=U(x)+iV(x)$ be a harmonic function such
that $F(\mathbb{B}^n)\subset\mathbb{D}$, $F(0)=a+bi$. Then, for
$0<r<1$ and $\o\in S$,
$$U(r\o)\le U_{a,b,r}(rN)$$
with equality at some point $r\o$ if and only if
$F(x)=F_{a,b,r}(xA)$, where $A$ is an orthogonal matrix such that
$r\o A=rN$, $U_{a,b,r}$ is defined as \eqref{Udy} and \eqref{UVtd},
$F_{a,b,r}$ is defined as \eqref{Fabrdy}. Further,
$U(x)<U_{a,b,r}(rN)$ for $|x|<r$. \el

\bp Step 1: First the case that $r\o=rN$ will be proved. Let
$0<\tilde{r}<1$ be fixed. Construct function
$$G(x)=F(\tilde{r}x)\quad\mbox{for}\quad x\in\ov{\mathbb{B}}^n.$$
$G(x)$ is harmonic on $\ov{\mathbb{B}}^n$ and $G(0)=a+bi$. Let
$G(x)=u(x)+iv(x)$. Then
\begin{equation} \label{Gtd}
\|u\|_{\infty}\le 1,\quad \int_S u(\o)d\s=a,\quad
\int_S\sqrt{1-u^2(\o)}d\s\ge\int_S |v(\o)|d\s\ge\lt|\int_S
v(\o)d\s\rt|=|b|.
\end{equation}
So by \eqref{uyq} we know that $u\in{\mc U}_{a,|b|}$ and by Lemma
\ref{lemmauabr}, we have $u(rN)\le U_{a,b,r}(rN)$ with equality if
and only if $u(\o)=u_{a,|b|,r}(\o)$ almost everywhere on $S$. For
$u_{a,|b|,r}(\o)$, by \eqref{u0td} and \eqref{b=0u0dy} we have
\begin{equation} \label{jfb}
\int_S \sqrt{1-u^2_{a,|b|,r}(\o)}d\s=|b|.
\end{equation}
If $u(\o)=u_{a,|b|,r}(\o)$ almost everywhere on $S$, then by
\eqref{Udy} and \eqref{UVtd}, we have
$$u(x)=U_{a,|b|,r}(x)=U_{a,b,r}(x)\ \ \mbox{for $x\in\mathbb{B}^n$};$$
and by \eqref{vdy}, we have
\begin{equation} \label{vbds}
v_{a,|b|,r}(\o)=\sqrt{1-u^2_{a,|b|,r}(\o)}=\sqrt{1-u^2(\o)}.
\end{equation}
Note that by \eqref{Gtd}, \eqref{jfb} and \eqref{vbds} we have
$$|b|=\int_S
v_{a,|b|,r}(\o)d\s\geq\int_S|v(\o)|d\s\geq\lt|\int_S
v(\o)d\s\rt|=|b|.$$ Then $$v(\o)=v_{a,|b|,r}(\o)\quad\mbox{almost
everywhere on}\quad S \quad\mbox{when}\quad b\geq0,$$
$$v(\o)=-v_{a,|b|,r}(\o)\quad\mbox{almost
everywhere on}\quad S \quad\mbox{when}\quad b<0.$$ So
$$v(x)=V_{a,b,r}(x)\ \ \mbox{for $x\in\mathbb{B}^n$}.$$

For $G(x)=u(x)+iv(x)$, it is proved that $u(rN)\le U_{a,b,r}(rN)$
with equality if and only if $G(x)=F_{a,b,r}(x)$. Now let
$\tilde{r}\rightarrow1$. Note that
$$\lim_{\tilde{r}\rightarrow1}G(x)=\lim_{\tilde{r}\rightarrow1}F(\tilde{r}x)=F(x),
\quad \lim_{\tilde{r}\rightarrow1}u(rN)=U(rN).$$ Then by the result
for $G(x)$, we have $U(rN)\le U_{a,b,r}(rN)$ with equality if and
only if $F(x)=F_{a,b,r}(x)$.

Step 2: Now we prove the case that $r\o\neq rN$. Construct function
$$\tilde{F}(x)=F(xA^{-1})\quad\mbox{for}\quad x\in\mathbb{B}^n,$$
where $A$ is an orthogonal matrix such that $r\o A=rN$ and $A^{-1}$
is the inverse matrix of $A$. By \cite{Axler}, we know that
$\tilde{F}(x)$ is also a harmonic function. Let
$\tilde{F}(x)=\tilde{U}(x)+i\tilde{V}(x)$. Note that
$\tilde{F}(0)=a+bi$. Then by the result of step 1, we have
$\tilde{U}(rN)\le U_{a,b,r}(rN)$ with equality if and only if
$\tilde{F}(x)=F_{a,b,r}(x)$. Note that
$\tilde{U}(rN)=U(rNA^{-1})=U(r\o)$ and $\tilde{F}(x)=F(xA^{-1})$.
Thus $U(r\o)\le U_{a,b,r}(rN)$ with equality if and only if
$F(xA^{-1})=F_{a,b,r}(x)$. It is just that $U(r\o)\le U_{a,b,r}(rN)$
with equality if and only if $F(x)=F_{a,b,r}(xA)$.

Step 3: We will show that $U(x)<U_{a,b,r}(rN)$ for $|x|<r$. By the
result of step 2 and the maximum principle, we have $U(x)\leq
U_{a,b,r}(rN)$ for $|x|\leq r$. If the equality holds for some $x_0$
with $|x_0|<r$, then $U(x)$ must be equal to $U_{a,b,r}(rN)$
identically for $|x|\leq r$. Note that if $U(rN)=U_{a,b,r}(rN)$,
then by the result of step 1, we have $U(x)=U_{a,b,r}(x)$. Thus
$U_{a,b,r}(x)\equiv U_{a,b,r}(rN)$ for $|x|\leq r$. However, it is
impossible since $U_{a,b,r}$ is not a constant. The proof of the
lemma is complete.
 \ep

\bl\label{lemma4} For fixed $0<r<1$ and $x\in \mathbb{B}^n$,
$F_{a,b,r}(x)$ is defined as \eqref{Fabrdy}. Then $F_{a,b,r}(x)$, as
a function of variables $a$ and $b$, is analytic in the real sense
on the open half disk $\{(a,b):b>0, a^2+b^2<1\}$ and is continuous
to the real diameter. \el

\bp Let $0<r<1$ and $x\in \mathbb{B}^n$ be fixed. It is obvious that
$F_{a,b,r}(x)$ is analytic in the real sense on the open half disk,
since it is determined there by the functions $\l(r,a,b)$ and
$\mu(r,a,b)$ formulated in Lemma 1, which are analytic in the real
sense on the open half disk $\{(a,b):b>0, a^2+b^2<1\}$. We only need
to prove that $F_{a,b,r}(x)$ is continuous at the points of the real
diameter. Note that \eqref{Fabrdy}. Then we only need to prove that
$U_{a,b,r}(x)$ and $V_{a,b,r}(x)$ are continuous at the points of
the real diameter.

Let $-1<a_0<1$ be given. We want to prove that $U_{a,b,r}(x)$ and
$V_{a,b,r}(x)$ is continuous at $(a_0,0)$. It is just to prove that
$U_{a,b,r}(x)\to U_{a_0,0,r}(x)$ and $V_{a,b,r}(x)\to
V_{a_0,0,r}(x)$ as $(a,b)\to (a_0,0)$.

Step 1: For the case that $(a,b)\to (a_0,0)$ with $b=0$, by
\eqref{Udy} and \eqref{Vdy}, we only need to prove $u_{a,0,r}(\o)\to
u_{a_0,0,r}(\o)$ almost everywhere on $S$ as $(a,0)\to (a_0,0)$.
Recall that
\begin{equation*}
u_{a,0,r}(\o)=
\begin{cases}
1, &\o\in S^+_{d_a};\\
0, &\o\in S_{d_a};\\
-1, &\o\in S^-_{d_a},
\end{cases}
\end{equation*}
where $S^+_{d_a}$, $S_{d_a}$ and $S^-_{d_a}$ are defined as
\eqref{Sd1}, \eqref{Sd2} and \eqref{Sd3}. This shows that
$u_{a,0,r}(\o)\to u_{a_0,0,r}(\o)$ almost everywhere on $S$ as
$(a,0)\to (a_0,0)$.

Step 2: For the case that $(a,b)\to (a_0,0)$ with $b>0$, by
\eqref{Udy} and \eqref{Vdy}, we only need to prove $u_{a,b,r}(\o)\to
u_{a_0,0,r}(\o)$ for any $\o\in S$ as $(a,b)\to (a_0,0)$ with $b>0$.

First we want to prove that $\mu(r,a,b)\to 0$ as $(a,b)\to (a_0,0)$
with $b>0$, where $\mu(r,a,b)$ is defined as $\mu(a,b)$ in
\eqref{muab}. Assume that $\mu(r,a,b)\nrightarrow 0$ as $(a,b)\to
(a_0,0)$ with $b>0$. Then there exists a sequence
$(a_k,b_k)\to(a_0,0)$ with $b_k>0$ such that $\mu_k=\mu(r,a_k,b_k)$
has a positive lower bound since $\mu(r,a,b)>0$. Then by
\eqref{muab} and \eqref{RIdy}, we have
$$\int_S\lt(1+\frac{1}{\mu^2_k}\lt(\frac{1}{|rN-\o|^n}-\l_k\rt)^2\rt)^{-1/2}d\s
=I(r,\l_k,\mu_k)=b_k\to0,$$ where $\l_k=\l(r,a_k,b_k)$, $\l(r,a,b)$
is defined as $\l(\mu(a,b),a)$ in \eqref{muab}. Thus
$\l_k\to\infty$. Assume that $\l_k\to+\iy$. Then by \eqref{u0dy},
\eqref{Ady} and \eqref{u0td}, we obtain that
$$u_{a_k,b_k,r}(\o)=\frac{\frac{1}{\mu_k}\lt(\frac{1}{|rN-\o|^n}-\l_k\rt)}
{\lt(1+\frac{1}{\mu^2_k}\lt(\frac{1}{|rN-\o|^n}-\l_k\rt)^2\rt)^{1/2}}\to-1,$$
uniformly for $\o\in S$, and  $a_k\to -1$, a contradiction.

Now, we want to prove that
$$\l(r,a,b)\to\l_0=\frac{1}{J^n_{a_0}}$$
as $(a,b)\to(a_0,0)$ with $b>0$, where $J^n_{a_0}$ is defined as
\eqref{xy4}. In contrary, assume that $\l(r,a,b)\nrightarrow\l_0$ as
$(a,b)\to(a_0,0)$ with $b>0$. Then there is a sequence
$(a_k,b_k)\to(a_0,0)$ with $b_k>0$ such that
$\l_k=\l(r,a_k,b_k)\to\l'\ne\l_0$. If $\l'=\iy$, then, as above,
$|a_k|\to1$, a contradiction. In the case that $\l'$ is finite, by
\eqref{u0dy}, \eqref{Ady} and \eqref{u0td} we have
$$u_{a_k,b_k,r}(\o)=\frac{\frac{1}{|rN-\o|^n}-\l_k}
{\lt(\mu^2_k\lt(\frac{1}{|rN-\o|^n}-\l_k\rt)^2\rt)^{1/2}}\to\mbox{sgn}\lt\{\frac1{|rN-\o|^n}
-\l'\rt\},$$
$$a_k=\int_S u_{a_k,b_k,r}(\o)d\s\to
\int_S \mbox{sgn}\lt\{\frac1{|rN-\o|^n} -\l'\rt\}d\s$$
\begin{equation*}
=
\begin{cases}
-1, &\l'\ge1/(1-r)^n;\\
1, &\l'\le1/(1+r)^n;\\
a', &\l'=1/J^n_{a'},\ -1<a'<1,\ a'\ne a_0.
\end{cases}
\end{equation*}
This contradicts $a_k\to a_0$.

It is proved that $\mu(r,a,b)\to0$ and $\l(r,a,b)\to\l_0$ as
$(a,b)\to(a_0,0)$ with $b>0$. Thus,
$$u_{a,b,r}(\o)\to\mbox{sgn}\lt\{\frac1{|rN-\o|^n}
-\l_0\rt\}=u_{a_0,0,r}(\o).$$

Step 3: For the case that $(a,b)\to (a_0,0)$ with $b<0$, by the
result of step 2, we know that $U_{a,-b,r}(x)\to U_{a_0,0,r}(x)$ and
$V_{a,-b,r}(x)\to V_{a_0,0,r}(x)$ as $(a,-b)\to (a_0,0)$ with
$-b>0$. Note that $U_{a,-b,r}(x)=U_{a,b,r}(x)$,
$V_{a,-b,r}(x)=-V_{a,b,r}(x)$ and $V_{a_0,0,r}(x)\equiv0$. Then we
have $U_{a,b,r}(x)\to U_{a_0,0,r}(x)$ and $V_{a,b,r}(x)\to
V_{a_0,0,r}(x)=0$ as $(a,b)\to(a_0,0)$ with $b<0$.

It is proved that $U_{a,b,r}(x)$ and $V_{a,b,r}(x)$ is continuous at
$(a_0,0)$. The lemma is proved.
 \ep

\section{Main results}

For $-\pi\le\b\le\pi$ and real number $\d$, denote the straight line
$l(\b,\d)$ and closed half plane $P(\b,\d)$ by
$$l(\b,\d)=\{w=u+iv: \mbox{Re}\{we^{-i\b}\}=u\cos\b+v\sin\b=\d\}$$
and
$$P(\b,\d)=\{w=u+iv: \mbox{Re}\{we^{-i\b}\}=u\cos\b+v\sin\b\leq\d\}.$$

\bt Let $0<r<1$ and $0\le\r<1$. Denote
$$P_\b=P(\b,U_{\r\cos\b,-\r\sin\b,r}(rN)),\quad
l_\b=l(\b,U_{\r\cos\b,-\r\sin\b,r}(rN))),$$ and define
\begin{equation}\label{Err}
E_{r,\r}=\bigcap_{-\pi\le\beta\le\pi}P_\b,
\end{equation}
$$
\G_{r,\r}=\{w:
w=f_{r,\r}(\b)=e^{i\b}F_{\r\cos\b,-\r\sin\b,r}(rN),-\pi\le\beta\le\pi\},$$
where $U_{\r\cos\b,-\r\sin\b,r}$ is defined as \eqref{Udy} and
\eqref{UVtd}, $F_{\r\cos\b,-\r\sin\b,r}$ is defined as
\eqref{Fabrdy}. \\
Then:
\begin{itemize}
\item[(1)] For any harmonic function $F$ such that $F(\mathbb{B}^n)\subset\mathbb{D}$ and $F(0)=\r$,
we have $F(\ov B_r)\su E_{r,\r}$;
\item[(2)] $E_{r,\r}$ is a closed convex domain
and symmetrical with respect to the real axis, and $\r$ is an
interior point of $E_{r,\r}$;
\item[(3)] $\Gamma_{r,\r}$ is a convex Jordan closed curve and $\p
E_{r,\r}=\Gamma_{r,\r}$;
\item[(4)] For any $w'\in E_{r,\r}$, there is a harmonic function
$F$ such that $F(\mathbb{B}^n)\subset\mathbb{D}$, $F(0)=\r$ and
$F(rN)=w'$.
\end{itemize}
\et

\bp (1) Denote
$$P'_\b=P(0,U_{\r\cos\b,-\r\sin\b,r}(rN)),\quad
l'_\b=l(0,U_{\r\cos\b,-\r\sin\b,r}(rN))).$$ $P_\b$ and $l_\b$ are
obtained from $P'_\b$ and $l'_\b$ by an anti-clockwise rotation of
angle $\b$.

Let $F$ be a harmonic function such that
$F(\mathbb{B}^n)\subset\mathbb{D}$ and $F(0)=\r$. For
$-\pi\le\beta\le\pi$, let $F_\b=e^{-i\b}F$. Then,
$F_\b(\mathbb{B}^n)\subset\mathbb{D}$ and
$F_\b(0)=\r(\cos\b-i\sin\b)$. Using lemma 3 to the harmonic function
$F_\b$, we have $F_\b(\ov B_r)\su P'_\b$ and, consequently, $F(\ov
B_r)\su P_\b$. This shows (1).

(2) It is obvious that $E_{r,\r}$ is a closed convex set and
symmetrical with respect to the real axis. We only need to prove
that $\r$ is an interior point of $E_{r,\r}$.

First we want to prove that $f_{r,\r}(\b)\in \p E_{r,\r}$ for
$-\pi\le\beta\le\pi$. $f_{r,\r}(\b)\in l_\b $ since
$F_{\r\cos\b,-\r\sin\b,r}(rN)\in l'_\b$. Let
$G(x)=e^{i\b}F_{\r\cos\b,-\r\sin\b,r}(x)$. The harmonic function $G$
satisfies the conditions $G(\mathbb{B}^n)\su\mathbb{D}$ and
$G(0)=\r$. By (1), $f_{r,\r}(\b)=G(rN)\in E_{r,\r}$. Note that
$E_{r,\r}\subset P_\b$, $l_\b=\p P_\b$ and $f_{r,\r}(\b)\in l_\b$
proved above. Then we have $f_{r,\r}(\b)\in \p E_{r,\r}$.

For $f_{r,\r}(0),f_{r,\r}(\pi),f_{r,\r}(\pi/2)$ and
$f_{r,\r}(-\pi/2)$, by lemma 3, we have
\begin{equation}\label{fro1}
f_{r,\r}(0)=F_{\r, 0,r}(rN)=U_{\r, 0,r}(rN)>U_{\r, 0,r}(0)=\r,
\end{equation}
\begin{equation}\label{fro2}
f_{r,\r}(\pi)=-F_{-\r, 0,r}(rN)=-U_{-\r, 0,r}(rN)<-U_{-\r,
0,r}(0)=\r,
\end{equation}
$$\mbox{Im}{f_{r,\r}(\pi/2)}=U_{0,-\r,r}(rN)=U_{0,\r,r}(rN)>U_{0,\r,r}(0)=0,$$
$$\mbox{Im}{f_{r,\r}(-\pi/2)}=-U_{0,\r,r}(rN)<-U_{0,\r,r}(0)=0.$$
Then $\r$ is an interior point of $E_{r,\r}$ since $E_{r,\r}$ is a
convex set.

(3) First we want to prove that $\G_{r,\r}$ is a Jordan closed
curve. $\G_{r,\r}$ is close and continuous by Lemma 4. Assume that
there exist $0<\b_1<\b_2<\pi$ such that
$w_0=f_{r,\r}(\b_1)=f_{r,\r}(\b_2)$. Then $\b_2-\b_1<\pi$ and $w_0$
is the vertex of the angular domain $P_{\b_1}\cap P_{\b_2}$.
Further, it is easy to see that $f_{r,\r}(\b)=w_0$ for
$\b_1<\b<\b_2$, since $l_\b\cap\p E_{r,\r}=w_0$ and $f_{r,\r}(\b)\in
l_\b\cap\p E_{r,\r}$. $f_{r,\r}(\b)$ is analytic on $(0,\pi)$ in the
real sense by Lemma 4. Then we have that $f_{r,\r}(\b)=w_0$ for
$0<\b<\pi$ and, by the continuity, $f_{r,\r}(0)=f_{r,\r}(\pi)=w_0$.
A contraction, since $f_{r,\r}(0)>f_{r,\r}(\pi)$ by \eqref{fro1} and
\eqref{fro2}. This shows that $\G^+_{r,\r}=\{w=f_{r,\r}(\b):
0\leq\b\leq\pi\}$ is a Jordan curve. By the same reason,
$\G^-_{r,\r}=\{w=f_{r,\r}(\b): -\pi\leq\b\leq0\}$ is also a Jordan
curve. Then $\G_{r,\r}$ is a Jordan closed curve.

For $-\pi\le\b\le\pi$, it is proved in (2) that $f_{r,\r}(\b)\in\p
E_{r,\r}$. Then $\G_{r,\r}\subset\p E_{r,\r}$. Note that $\p
E_{r,\r}$ must be a convex Jordan closed curve. Thus $\p
E_{r,\r}=\G_{r,\r}$.

(4)For $w'\in E_{r,\r}$, draw a straight line $l$ passing through
$w'$ and intersect $\p E_{r,\r}$ at $w_1$ and $w_2$. Let
$w'=k_1w_1+k_2w_2$ with $k_1,k_2\geq0$ and $k_1+k_2=1$. There are
two real numbers $\b_1$ and $\b_2$ such that $f_{r,\r}(\b_1)=w_1$
and $f_{r,\r}(\b_2)=w_2$. Then the harmonic function
$F=k_1e^{i\b_1}F_{\r\cos\b_1,-\r\sin\b_1,r}+k_2e^{i\b_2}F_{\r\cos\b_2,-\r\sin\b_2,r}$
satisfies $F(\mathbb{B}^n)\subset\mathbb{D}$, $F(0)=\r$ and
$F(rN)=w'$. The theorem is proved. \ep

When $\r=0$, we have a corollary as follows, which is coincident
with \eqref{HSlemma4}.

\bc Let $0<r<1$. For any harmonic mapping $F$ such that
$F(\mathbb{B}^n)\subset\mathbb{D}$ and $F(0)=0$, we have
$$F(\ov B_r)\su \ov D_{U(rN)},$$
where $U$ is the Poisson integral of the function that equals 1 on
$S^+$ and -1 on $S^-$.\ec

\bp By Theorem 2, we only need to prove that $E_{r,0}=\ov
D_{U(rN)}$. Further, by the definition of $E_{r,\r}$ in Theorem 2,
we only need to prove that $U_{0,0,r}(rN)=U(rN)$. Note that by
\eqref{b=0u0dy},
\begin{equation*}
u_{0,0,r}(\o)=
\begin{cases}
1, &\o\in S^+;\\
0, &\o\in S;\\
-1, &\o\in S^-.
\end{cases}
\end{equation*}
Then by \eqref{Udy} we know that $U_{0,0,r}(rN)=U(rN)$. The
corollary is proved. \ep

From Theorem 2, we obtain Theorem 1, which is the general version of
the above Theorem 2.

\newpage

\end{document}